\documentclass [11pt]{article}
\usepackage{amsthm}
\usepackage{amssymb}
\usepackage{epsfig}
\usepackage{srcltx}
\usepackage[dvipsnames,usenames]{color}

\newcounter{contador}
\newtheorem{propo}[contador]{Proposition}
\newtheorem{teo}[contador]{Theorem}
\newtheorem{lem}[contador]{Lemma}
\newtheorem{nota}[contador]{Remark}

\newtheorem{corol}[contador]{Corollary}

\newcommand{\rec}{\noindent}    
\renewcommand{\proof}{\noindent \textsl{Proof.} }  
\renewcommand{\qed}{\, \hfill\rule[-1mm]{2mm}{3.2mm}} 

\newcommand{\R}{{\mathbb R}}
\newcommand{\C}{{\mathbb C}}
\newcommand{\N}{{\mathbb N}}

\title{On periodic solutions of $2$--periodic Lyness'
equations\footnote{The second author is partially supported by
Spain's Ministry of Science and Technology (MCYT) through grant
DPI2011-25822. CoDALab group is supported by the Catalonia's
government through the SGR program. The support of DMA3's Terrassa
Campus Section is also acknowledged.}}

\author{Guy Bastien$^{(1,4)}$, V\'{\i}ctor Ma\~{n}osa$^{(2)}$ and  Marc Rogalski$^{(3,4)}$
  \\*[.1truecm]
{\small \textsl{$^{(1)}$ Institut Math\'ematique de Jussieu,}}
\\*[-.25truecm] {\small \textsl{Universit\'e Paris 6 and CNRS, France.}}
\\*[-.25truecm] {\small \textsl{bastien@math.jussieu.fr}}
\\*[-.25truecm]
\\*[-.25truecm] {\small \textsl{$^{(2)}$ Departament de Matem\`{a}tica Aplicada III,}}
\\*[-.25truecm] {\small \textsl{Control, Dynamics and Applications Group}}
\\*[-.25truecm] {\small \textsl{Universitat Polit\`{e}cnica de Catalunya}}
\\*[-.25truecm] {\small \textsl{Colom 1, 08222 Terrassa, Spain}}
\\*[-.25truecm] {\small \textsl{victor.manosa@upc.edu}}
\\*[-.25truecm]
\\*[-.25truecm] {\small \textsl{$^{(3)}$ Laboratoire Paul Painlev\'e, }}
\\*[-.25truecm] {\small \textsl{Universit\'e de Lille 1 and CNRS, France.}}
\\*[-.25truecm] {\small \textsl{marc.rogalski@upmc.fr}}
\\*[-.25truecm]
\\*[-.25truecm] {\small \textsl{$^{(4)}$ Universit\'e Paris 6 and CNRS, }}
\\*[-.25truecm] {\small \textsl{4 pl. Jussieu, 75005 Paris,
France.}}}

\begin{document}
\maketitle

\begin{abstract}
We study the existence of periodic solutions of the
\emph{non--autonomous periodic Lyness' recurrence}
$u_{n+2}=(a_n+u_{n+1})/u_n,$ where $\{a_n\}_n$ is a cycle with
positive values $a$,$b$ and with positive initial conditions. It is
known that for $a=b=1$  all the sequences generated by this
recurrence are $5$--periodic.  We prove that for each pair
$(a,b)\neq(1,1)$ there are infinitely many initial conditions giving
rise to periodic sequences, and that the family of recurrences have
almost all the even periods. If $a\neq b$, then any odd period,
except $1$, appears.
\end{abstract}

\rec {\sl 2000 Mathematics Subject Classification:}
\texttt{39A20,39A11} primary;  \texttt{14H52,14H70} secondary.

\rec {\sl Keywords:} Difference equations with periodic
coefficients; elliptic curves; Lyness' type equations;  QRT maps;
rotation number; periodic orbits.

\rec {\sl Running title:} $2$-periodic Lyness equations

\newpage

\section{Introduction and main results}

The study of periodic orbits in periodic non--autonomous discrete
dynamical systems is a classical topic that, mainly driven by some
conjectures in mathematical biology, has attracted again the
researcher's interest in the last years. This is because these ones
are good models for describing the dynamics of biological systems
under periodic fluctuations whether due to external disturbances or
effects of seasonality  (see \cite{BHS,CH,CH2,ES1,ES2,Sv} for
instance and references therein). On the other hand, the existence
of discrete integrable systems with periodic coefficients is
currently being the focus of some mathematical physics research both
on those systems belonging to the celebrated \emph{QRT family of
maps} (\cite{D}), or on other families as the HKY one (see
\cite{GRTW,RGW} and \cite{GRT2011}, respectively). Our work is in
the interface of these two streams, since the main goal of the paper
is to characterize the global periodic structure of a particular
family of QRT maps: the one associated to the $2$-periodic Lyness
equations.

\medskip

Nowadays, the dynamics of the \emph{autonomous Lyness' difference
equation}
\begin{equation}\label{Lynesseq}
u_{n+2}=\displaystyle{\frac{a+u_{n+1}}{u_{n}}} \mbox{ with } a>0,
\mbox{ and } u_1,u_2>0.
\end{equation}
is completely understood after the research done in \cite{BR,BC} and
\cite{Z} (see also \cite{GMX11}). In summary, the dynamics of
equation (\ref{Lynesseq}) can be studied through the dynamics of the
map $F_{a}(x,y)=\left(y,(a+y)/{x}\right)$. This map has a first
integral $V_a$ such that their level sets in
$\mathcal{Q}^+:=\{(x,y)\in\R^2:x>0,y>0\}$, except the one
corresponding to the unique fixed point, are the ovals of some
\emph{elliptic curves} of the form
\begin{equation}\label{elipautonom}
\{V_a(x,y)=h\}=\{(x+1)(y+1)(x+y+a)-hxy=0\}.
\end{equation}
The action of the map on each of the above curves can be described
in terms of the \emph{group law} of them: in fact, the Lyness' one
is a particular case of the well known family of QRT maps, \cite{D}.
In particular all possible periods of the recurrences generated by
(\ref{Lynesseq}) are known, and for any $a\notin\{0,1\}$ infinitely
many different prime periods appear (the cases $a=0$ and $1$ are
globally periodic with periods $6$ and $5$ respectively).

Recently, there has been some progress concerning the study of the
\emph{non--autonomous} periodic Lyness' equations
\begin{equation}\label{eq}
u_{n+2}\,=\,\frac{a_n+u_{n+1}}{u_n},
\end{equation}
when $\{a_n\}_n$ is a $k$-periodic sequence taking positive values,
 and the initial conditions $u_1,u_2$ are, as well, positive (see
 \cite{CGM11a,CGM11b,CZ,JKN,KN}, and also \cite{GRTW,RGW}). These works focus on some
 qualitative aspects of the dynamics like
 persistence, stability,  as well as some integrability issues.
In this paper we will focus on the characterization of periodic
solutions of the $2$--periodic case given by
\begin{equation}\label{k=2}a_n\,=\,\left\{\begin{array}{lllr}
a&\mbox{for}&n=2\ell+1,\\ b&\mbox{for}&\,n=2\ell,
\end{array}\right.\end{equation}
where $\ell\in\mathbb{N}$ and $a>0,b>0$.

The solutions of equation (\ref{eq}) can be studied through the
dynamics given by the \emph{{composition map}}
\begin{equation}\label{FBI}F_{b,a}(x,y):=(F_{b}\circ F_{a})(x,y)=\Big({a+y\over x},{a+bx+y\over xy}\Big),\end{equation}
since the relation between the terms of recurrence
(\ref{eq})--(\ref{k=2}) and the iterates of the composition map is
given by
\begin{equation}\label{Fbaab}
(u_{2n+1},u_{2n+2})=F_{b,a}(u_{2n-1},u_{2n}), \mbox{ and }
(u_{2n+2},u_{2n+3})=F_{a,b}(u_{2n},u_{2n+1}),
\end{equation}
where  $(u_1,u_2)\in{\mathcal{Q}^+}$ and $n\ge 1$.

The map $F_{b,a}$ is \emph{integrable}  in the sense that it has a
first integral, \cite{JKN}. Here we use the following one:
\begin{equation}\label{invariantVba}
V_{b,a}(x,y)=\frac{(bx+a)(ay+b)(ax+by+ab)}{xy}.
\end{equation}
This means that each map $F_{b,a}$ preserves a foliation of the
plane given by cubic curves which are, generically, elliptic. The
composition maps $F_{b,a}$ associated to the $2$--periodic Lyness'
equations are also particular cases of QRT maps. It is interesting,
however, to notice that the $2$--periodic case is one of the few
\emph{rationally integrable}  or even \emph{meromorphically
integrable} $k$--periodic Lyness' equations (see \cite{CGM11b,CZ}).
The maps $F_{b,a}$ are also a generalization of the map considered
in \cite{ER}.

It is known that a map $F_{b,a}$ has a unique fixed point
$(x_c,y_c)\in \mathcal{Q}^+$ given by the solution of the system
\begin{equation}\label{puntfix}
 \left\{
   \begin{array}{l}
     x^2=a+y, \\
     y^2=b+x.
   \end{array}
 \right.
\end{equation}
which corresponds to the unique global minimum of $V_{b,a}$ in
${\mathcal{Q}^+}$. Furthermore, setting $h_c:=\{V_{b,a}(x_c,y_c)\}$,
the level sets
$\mathcal{C}_h^+:=\left\{\{V_{b,a}=h\}\cap\mathcal{Q}^+\right.$ for
$\left. h>h_c\right\}$ are closed curves and the dynamics of
$F_{b,a}$ restricted to these sets is \emph{conjugate to a rotation
on the unit circle} with associated \emph{rotation number}
$\theta_{b,a}(h)$, \cite{CGM11a}. In Section \ref{secrotacio} we
give an alternative proof of this fact (see Corollary
\ref{theyarerotations}). In this paper we study some properties of
the rotation number in order to characterize the possible periods
that can appear in the family of recurrences~(\ref{eq}).

\medskip

The main results of the paper are the following.

\newpage

  \begin{teo}\label{periodsinI} Set
$$I(a,b):=\left\langle \sigma(a,b),{2\over5} \right\rangle,$$
where $\langle c,d\rangle=(\min(c,d),\max(c,d))$, and $$ \sigma(a,b)
:=\frac{1}{2\pi}\arccos\left(\frac{1}{2}\left[-1-\displaystyle{\frac{a+b\,x_c}{x_c\,(b+x_c)}}\right]\right)
=\displaystyle{ \frac{1}{2\pi}}\arccos\left(\displaystyle{
\frac{1}{2}}\left[-2+\displaystyle{\frac{1}{x_c\,y_c}}\right]\right).
$$
For any fixed $a,b>0$, and any value $\theta\in I(a,b)$, there exist
at least an oval of the form $\mathcal{C}_h^+$ such that the map
$F_{b,a}$ restricted to the this oval  is conjugate to a rotation,
with a rotation number $\theta_{b,a}(h)=\theta$.
\end{teo}

\vspace{0.2cm}

 Notice that $F_{1,1}$ is the doubling of the well--known globally
autonomous Lyness' map $F_1$ which is globally $5$--periodic and
whose rotation number function takes the constant value
$1/5$\footnote{Here we adopt the notation in \cite{BR} and \cite{Z}
concerning the determination of the rotation number.  Notice that in
\cite{CGM11a} it is adopted the determination given by
$1-\theta_{b,a}(h)$.}, and therefore $I(1,1)=\{2/5\}$. Notice also,
that there are other values of $a$ and $b$ for which
$I(a,b)=\{2/5\}$, as can be seen in Proposition
\ref{propocondiciosigma}.

\begin{teo}\label{teo2} Consider the family of maps $F_{b,a}$ given
in (\ref{FBI}) for $a,b>0$.
\begin{itemize}
\item[(i)] If  $(a,b)\neq(1,1)$, then there exists a  value
$p_0(a,b)\in\N$ such that for any $p>p_0(a,b)$ there exist at least
a continuum of initial conditions in ${{\cal Q}^+}$ (an oval
$\mathcal{C}_h^+$) giving rise to $p$--periodic orbits of $F_{b,a}.$
Moreover, when $\sigma(a,b)\neq 2/5$ the number $p_0(a,b)$ is
computable.

  \item[(ii)] For each number $\theta$ in $(1/3,1/2)$ there exists some $a>0$
and $b>0$ and at least an oval $\mathcal{C}_h^+$, such that the
action of $F_{b,a}$ restricted to this oval is conjugate to a
rotation with rotation number $\theta_{b,a}(h)=\theta.$ In
particular, for all the irreducible rational numbers $q/p\in
(1/3,1/2)$, there exist periodic orbits of $F_{b,a}$ of prime period
$p.$
  \item[(iii)] The set of periods arising in the family
$\{F_{b,a},\,a>0,b>0\}$ restricted to ${\mathcal{Q}^+}$ contains all
prime periods except $2$, $3$, $4$, $6$ and $10$.
\end{itemize}
\end{teo}

In fact, as we will see in Section \ref{provamainresults}, the prime
periods $2$ and $3$ do not appear for any $a$ and $b$  in the whole
domain of definition of the dynamical system defined by $F_{b,a}$,
but periods $4$, $6$, and $10$ appear for some $a,b>0$ and some
initial conditions in ${\cal G}\setminus{\cal Q}^+$.

\begin{corol}\label{reclyness}
Consider the $2$--periodic  Lyness' recurrence
(\ref{eq})-(\ref{k=2}) for $a>0$, $b>0$ and positive initial
conditions $u_1$ and $u_2$.
\begin{itemize}
\item[(i)] If  $(a,b)\neq(1,1)$, then there exists a value $p_0(a,b)\in\N$
 such that
for any $p>p_0(a,b)$ there exist continua of initial conditions
giving rise to $2p$--periodic sequences. Moreover, when
$\sigma(a,b)\neq 2/5$ the number $p_0(a,b)$ is computable.

\item[(ii)]
The set of prime periods arising when $(a,b)\in(0,\infty)^2$ and
positive initial conditions are considered contains all the even
numbers except $4$, $6$, $8$,  $12$ and $20$. If $a\neq b$, then it
does not appear any odd period, except $1$.
\end{itemize}
\end{corol}

 With our tools, the number $p_0(a,b)$ in the statement
(i) of both Theorem \ref{teo2} and Corollary \ref{reclyness} is only
computable if $\sigma(a,b)\neq 2/5$. As can be seen from Proposition
\ref{propocondiciosigma} or Corollary \ref{corolgamma}, this happens
in an open an dense subset of the parameter space, and therefore
$p_0(a,b)$ is a \emph{generically} computable value.

Observe that Theorem \ref{teo2}  characterizes all the prime periods
that can appear for the iterations of the maps $F_{b,a}$ in
${\mathcal{Q}}^+$. In fact it characterizes the set of periods of
any planar system of first order difference equations associated to
any map conjugate to $F_{b,a}$, like the two families of first order
systems of difference equations
$$\left\{
    \begin{array}{ll}
      u_{n+1}u_n=a+v_n, \\
     v_{n+1}v_n=b+u_{n+1},
    \end{array}\right.\,\mbox{ and }\, \left\{
    \begin{array}{ll}
      u_{n+1}u_n=\alpha(1+v_n), \\
     v_{n+1}v_n=\beta(1+u_{n+1}).
    \end{array}\right.
$$ for some choices of
$a,b,\alpha,\beta>0$, and $u_1>0$, $v_1>0$. The second one  is
associated to the map $G$, defined in (\ref{funcioG}) below.

The paper is structured as follows. In Section \ref{secrotacio} we
see that on the level sets  $\mathcal{C}_h^+$, the action of the map
$F_{b,a}$ can be seen as a linear action of a birational map on an
elliptic curve, and thus we reobtain that on these level sets it is
conjugate to a rotation (these also follows from the fact of being a
QRT map, but the notation and elements introduced here will be
useful in our further analysis of the rotation number function). The
main part of this section is devoted to proof the elliptic nature of
the invariant curves on $\mathbb{R}^2$,
$\mathcal{C}_h:=\{V_{b,a}=h\}$, as well as to derive a
\emph{Weierstrass normal form} representation of both the invariant
curves and the map. This one, will be a key step in our approach to
the asymptotic behavior of the rotation number function when $h$
tends to infinity. This is done in Section \ref{behav}, which is
devoted to study the limit of the rotation function. In particular,
we obtain that the asymptotic behavior of this function at the
energy level corresponding to the fixed point and at the infinity do
not coincide. It is worth noticing that further tools to study the
rotation intervals have to be developed to face the problem of
finding the set of periods when these limits coincide. The main
results, as well as some other ones concerning periodic orbits are
proved in Section~\ref{provamainresults}.

\section{$F_{b,a}$ as a linear action in terms of
the group law of a cubic}\label{secrotacio}
\subsection{An overview from an algebraic geometric
viewpoint}\label{overview}

In this section we will see $F_{b,a}$ as a linear action in terms of
the group law of the cubic. The level sets $\{V_{b,a}=h\}$ are given
by the cubic curves in $\R^2$:
$$
{\cal C}_h=\{(bx+a)(ay+b)(ax+by+ab)-hxy=0\}.
$$
These curves ${\cal C}_h$, in homogeneous coordinates $[x:y:t]\in\C
P^2$, write as
$$
\widetilde{\cal C}_h=\{(bx+at)(ay+bt)(ax+by+abt)-hxyt=0\}.
$$
Observe that there are three infinite points at infinity which are
common to all the above curves
$$H=[1:0:0];\,\,\,V=[0:1:0];\,\,\,D=[b:-a:0].$$
Notice that none of these points is an inflection point.

The map $F_{b,a}$ extends naturally to $\C P^2$ as
$$
{\widetilde{F}}_{b,a}\left([x:y:t]\right)=\left[ayt+y^2:at^2+bxt+yt:xy\right],
$$ and leaves the curves $\widetilde{\cal C}_h$ invariant.

Notice  that  only a finite number curves $\widetilde{\cal C}_h$ are
singular, and thus correspond to non--elliptic curves ${\cal C}_h$,
since the discriminant is a polynomial in $h$ (see Section
\ref{elmodif}). In particular, we will prove that for all the energy
levels $h>h_c$, the curves $\widetilde{\cal C}_h$ are non-singular.

\begin{propo}\label{iselliptic1}
If $a>0$ and $b>0$, and for all $h>h_c$, the curves $\widetilde{\cal
C}_h$ are elliptic.
\end{propo}

The above result will be a direct consequence of Proposition
\ref{iselliptic} (proved in Section \ref{provaellipt}), and it
implies that in ${\cal Q}^+\setminus\{(x_c,y_c)\}$ the map $F_{b,a}$
is a \emph{birational transformation on an elliptic curve}, and
therefore it can be expressed as a linear action in terms of the
group law of the curve (\cite[Theorem 3]{JRV}). Indeed, for those
cases such that $\widetilde{\cal C}_h$ is elliptic, if we consider
$P_0=[x_0:y_0:1]\in\widetilde{\cal C}_h$, and taking the infinite
point $V$ as the zero element of $\widetilde{\cal C}_h$, then
$P_0+H$ can be computed as follows: take the horizontal line passing
through $P_0$ and $H$, it cuts $\widetilde{\cal C}_h$ in a second
point $(P_0*H)=(x_1,y_0)$ (where $M*N$ denotes the point where the
line $\overline{MN}$ cuts again the cubic). A computation shows that
$x_1=(a+y_0)/x_0$. The vertical line passing through $P_0*H$ cuts
$\widetilde{\cal C}_h$ at a new point $P_1=[x_1:y_1:1]=P_0+H$. Again
a computation shows that $y_1=(a+bx_0+y_0)/(x_0y_0)$. Hence
$(x_1,y_1)=F_{b,a}(x_0,y_0)$, and we get the following result:

\begin{propo}\label{lemasuma}
For each value of $h$ such that $\widetilde{\cal C}_h$ is an
elliptic curve, ${\widetilde{F}}_{{b,a}{|_{\widetilde{\cal
C}_h}}}(P)=P\,+\,H$, where $+$ is the addition of the group law of
 $\widetilde{\cal C}_h$ taking the infinite point $V$ as the zero
element.
\end{propo}

So we have the following immediate corollary

\begin{corol}\label{theyarerotations}
On each level set $\mathcal{C}_h^+$
 the
map $F_{b,a}$ is conjugate to a rotation.
\end{corol}

\begin{nota}\label{notaa1} Notice that the zero element $V$ of the group law is not an inflection
point of the cubic. So the relation of collinearity for three points
$A,B,C$ does not mean that $A+B+C=V$, but $A+B+C=V*V$ (an inflection
point $U$ satisfies $U*U=U$). Another useful relation is
$-P=P*(V*V)$, where $V*V=R:=[-a/b:0:1]$.
\end{nota}

The above results are a consequence of the ones proved in the next
section. Now we are able to give a first result about admissible
periodic orbits of $F_{b,a}$.

\begin{propo}\label{finit}
For any admissible period $p$ of $F_{b,a}$, there is only a finite
number of level curves of the admissible invariant levels sets
$\mathcal{C}_h$, which may have this period. \end{propo}

\proof As we have noticed, there is only a finite number of energy
levels $\{V_{b,a}=h\}$, which correspond to non--elliptic curves.
For the rest of values (those such that $\widetilde{\cal C}_h$ is an
elliptic curve) we have the following: if there is a periodic orbit
with minimal period $2n$ in the curve
$\widetilde{\mathcal{C}}_h\cap{\cal Q}^+$ with $h>h_c$, since
$F_{b,a}$ is conjugate to a rotation, then
$\widetilde{\mathcal{C}}_h$  must be filled by with periodic orbits
of minimal period $2n$, and also, taking into account the natural
extension $\widetilde{F}_{b,a}$ (defined in Proposition
\ref{lemasuma}) we have $2nH=V$ or, in other words, if and only if
$nH=-nH$.  Using Remark \ref{notaa1} we have that $-nH=nH*R$, where
$R=(-a/b,0,1)$. So there exists $2n$--periodic orbits if and only if
$nH=nH*R$, that is, if and only if $R$ is on the tangent line to the
curve $\widetilde{{\cal C}}_h$ at the point $nH$.  Computing this
tangent line and imposing this condition, we obtain that the energy
level of corresponding to this curve must be a zero of a polynomial
$P_{2n}(h)$, so it has a finite number of solutions.

Suppose that ${\cal C}_h$ is filled by periodic orbits of period
$p=2n+1$, then
$$(2n+1)H=V\Leftrightarrow nH+H=-nH \Leftrightarrow \left(nH*H\right)*V=nH*R,$$
which is again an algebraic relation on $h$ giving rise to a finite
number of solutions.\qed

\subsection{Proof of the ellipticity of the curves $\widetilde{\cal C}_h$
for $h>h_c$ and a Weierstrass normal form representation for
$F_{b,a}$}\label{provaellipt}

In this section we prove Proposition \ref{iselliptic1}. In order to
simplify the computations we introduce a preliminary change of
variables and parameters. It is easy to check from
(\ref{invariantVba}) that the integral $V_{b,a}$ admits the
following form
$${V_{b,a}(x,y)}/{(a^2b^2)}=\Big(\frac{x}{b}+\frac{a}{b^2}\Big)\Big(\frac{y}{a}+\frac{b}{a^2}\Big)\Big(1+\frac{x}{b}+\frac{y}{a}\Big)/\Big[\left(\frac{x}{b}\right)\left(\frac{y}{a}\right)\Big].$$

This means that  taking the new variables
\begin{equation}\label{changeparameters}
\left\{X:=\frac{x}{b},\,\,Y:=\frac{y}{a},\,\,\alpha:=\frac{a}{b^2},\,\,\beta:=\frac{b}{a^2},\right.
\end{equation}
the first integral is now given by
\begin{equation}\label{W}W(X,Y)=\frac{(X+\alpha)(Y+\beta)(1+X+Y)}{XY};\end{equation}
the map $F_{b,a}$ becomes
\begin{equation}\label{funcioG}
G(X,Y)=\Big({\alpha(1+Y)\over X},{\beta(X+\alpha Y+\alpha)\over
XY}\Big);\end{equation} and the curves $\widetilde{\cal C}_h$ are
brought into
$$
{\cal D}_L=\{(X+\alpha T)(Y+\beta T)(T+X+Y)-LXYT=0\},
$$
where $L=h/(a^2b^2)$. Observe that all the curves ${\cal D}_L$ have
three common infinite points given by $H=[1:0:0]$; $V=[0:1:0];$ and
$D=[1:-1:0],$ and again, none of these points is an inflection
point. The map $G$ extends to $\C P^2$ as
$$
\widetilde{G}\left([X:Y:T]\right)=\left[\alpha(YT+Y^2):\beta(XT+\alpha
YT+\alpha T^2):XY\right]
$$

The level $L_c=h_c/(a^2b^2)$ corresponds to the level of the global
minimum of the invariant $W$ in ${\cal Q}^+$. Proposition
\ref{iselliptic1} is, then, a direct consequence of the following
result, that will be proved in  Subsection \ref{elmodif}:

\begin{propo}\label{iselliptic}
If $\alpha>0$ and $\beta>0$, for all  $L>L_c=h_c/(a^2b^2)$, the
curves ${\cal D}_L$ are elliptic.
\end{propo}

Reasoning as in the case of $F_{b,a}$  we obtain the following
result.

\begin{propo}\label{l9}
For each $L$ such that ${\cal D}_L$ is an elliptic curve,
$\widetilde{G}_{|_{{\cal D}_L}}(P)=P\,+\,H$, where $+$ is the
addition of the group law of ${\cal D}_L$ taking the infinite point
$V$ as the zero element\footnote{ Notice that, again, the zero
element $V$ of the group law is not an inflection point of the cubic
(see Remark \ref{notaa1}).}.
\end{propo}

As a consequence of the Propositions \ref{iselliptic} and \ref{l9}
we obtain that the action of $G$ on each level ${\cal D}_L\cap {\cal
Q}^+$, with $L>L_c$ is conjugate to a rotation.

In order to prove Proposition \ref{iselliptic} as well as to
significatively simplify some of the proofs of the results
concerning   the behavior of $\theta_{b,a}(h)$ done in Section
\ref{behav},  we first present a new conjugation between the linear
action of $\widetilde{G}$ on ${\cal D}_L$ (and therefore for
$\widetilde{F}_{b,a}$ on $\widetilde{\cal C}_h$) and a linear action
on a standard Weierstrass normal form of the invariant curves. We
prove:

\begin{propo}\label{propo1}
\begin{itemize}
  \item[(i)] For any fixed $L$, the curves ${\cal D}_L$ have
  an associated Weierstass form given by
\begin{equation}\label{Wei}
\mathcal{E}_L=\{[x:y:t],\,y^2t=4\,x^3-g_2\,xt^2-g_3\,t^3\},
\end{equation} being
$$g_2=\frac{1}{192}\left({L}^{8}+\sum_{i=4}^{7}p_{i}(\alpha,\beta)L^i\right)\,
\mbox{ and }\,
g_3=\frac{1}{13824}\left(-L^{12}+\sum_{i=6}^{11}q_{i}(\alpha,\beta)L^i\right),$$
where
$$
\begin{array}{rl}
p_7(a,b)=&-4\left(\alpha+\beta+1
\right),  \\
p_6(a,b)=&2\left(3(\alpha-\beta)^2+2(\alpha+\beta)+3 \right),\\
p_5(a,b)=&-4 \left(\alpha+\beta-1 \right)  \left( {\alpha}^{2}-4
\beta\alpha+{\beta}^{2}-1 \right),\\
p_4(a,b)=&\left(\alpha+\beta-1
 \right) ^{4}.
\end{array}$$

\noindent and
$$
\begin{array}{rl}
q_{11}(a,b)=&6\left(\alpha+\beta+1
\right),\\
q_{10}(a,b)=&3\left(-5{\alpha}^{2}+2\alpha\beta-5{\beta}^{2}-6\alpha-6\beta-5\right)\\
q_{9}(a,b)=&4\left(5{\alpha}^{3}-12{\alpha}^{2}\beta-12\alpha{\beta}^{2}+5{
\beta}^{3}+3{\alpha}^{2}-3\alpha\beta+3{\beta}^{2}+3\alpha+3 \beta+5
\right)\\
q_{8}(a,b)=&3\left(
-5{\alpha}^{4}+16{\alpha}^{3}\beta-30{\alpha}^{2}{\beta}^{2}+16
\alpha{\beta}^{3}-5{\beta}^{4}+4{\alpha}^{3}\right.\\
&\left.-12{\alpha}^{2}
\beta-12\alpha{\beta}^{2}+4{\beta}^{3}+2{\alpha}^{2}-8\alpha
\beta+2{\beta}^{2}+4\alpha+4\beta-5\right)
\\
q_{7}(a,b)=&6 \left( {\alpha}^{2}-4\alpha\beta+{\beta}^{2}-1 \right)
 \left(\alpha+\beta-1 \right) ^{3}
\\
q_{6}(a,b)=&- \left(\alpha+\beta-1 \right) ^{6}
\end{array}$$
  \item[(ii)] For any $L>L_c$,  the curve $\mathcal{E}_L$ is an elliptic
  curve.
  \item[(iii)] For any $L>L_c$,  $\widetilde{G}_{|{\cal D}_L}$ is conjugate to the linear action
\begin{equation}\label{accio}
\widehat{G}_{|_{{\cal E}_L}}:P\mapsto P\,\star\,\widehat{H},
\end{equation}
where
\begin{equation}\label{eqhtilde}\widehat{H}=\left[\frac{1}{48}\left({L}^{2}-2\left(\alpha+\beta+1
\right) L+ \left( \alpha+\beta-1
 \right) ^{2}
\right)L^2\,:\,-\frac{1}{8}\alpha\beta
L^4\,:\,1\right]\end{equation} and $\star$ denotes the sum operation
of the cubic $\mathcal{E}_L$, taking the infinity point
$\widehat{V}=[0:1:0]$ as the zero of the group law.
\end{itemize}
\end{propo}

\vspace{0.2cm}

In the next subsections we give the proof of the above result.

\subsubsection{Computing the Weierstrass normal form of ${\cal
D}_L$}\label{weiers}

In this section we explicit the transformations which pass ${\cal
D}_L$ into their Weierstrass normal form. Although it is a standard
calculation, we prefer to include it because in the proof of
Proposition \ref{propo1}, we need the explicit expressions of some
of the intermediate transformations below.

\smallskip

 \noindent  \textsl{First transformation:} Following \cite[Section
I.3]{ST}, we take the following reference projective system: Let
$t=0$  be the tangent line to ${\cal D}_L$ at $V$. This tangent line
intersects  ${\cal D}_L$ at $R=[-\alpha:0:1]$, so  we will take the
reference $x=0$  to be the tangent line
 at $R$. Finally $y=0$ will be the $Y$ axis. The new
coordinates $x,y,t$ are, then, obtained via
\begin{equation}\label{firsttrans}
  x=\beta(1-\alpha)X+\alpha L Y+\alpha\beta(1-\alpha)T;\, y=X;
       t=X+\alpha T,
\end{equation}
 and the equation of the cubic becomes
\begin{equation}\label{firsteq}xy^2+y(a_1xt+a_2t^2)=a_3x^2t+a_4xt^2+a_5t^3,\end{equation}
where \begin{equation}\label{parama}\left\{
    \begin{array}{l}
      a_1=(\alpha-\beta-1-L)/L,\quad
      a_2=\beta\left((1-a)^2+\beta(1-\alpha)-L\right)/L,\quad
      a_3=-1/L^2,\\
      a_4=\left(2\beta(1-\alpha)-L(1+\beta)\right)/L^2,\quad
      a_5=\beta(\beta(1-\alpha)-L)(L-1+\alpha)/L^2.
    \end{array}\right.
\end{equation}

\noindent  \textsl{Second transformation:}  We multiply both members
of (\ref{firsteq}) by $x$, and we take the new  variables (denoted
again $X,\,Y,\,T$):
\begin{equation}\label{sectrans}
X=x;
       Y=xy/t;
       T=t,
\end{equation}
obtaining the new equation
\begin{equation}\label{seceq}T^2Y^2+(a_1X+a_2T)YT^2=a_3X^3T+a_4X^2T^2+a_5XT^3.\end{equation}

\noindent \textsl{Third transformation:} Removing the factor $T$ in
(\ref{seceq}) and taking the new variables $x,y,t$
\begin{equation}\label{thirdtrans}
x=X;\,y=Y+(a_1X+a_2T)/2;\, t=T,
\end{equation}
we obtain the new equation
$$ty^2=\lambda x^3+\mu x^2t+\nu xt^2+\gamma t^3,$$
where
\begin{equation}\label{parammu}\lambda=a_3,
\mu=a_4+{a_1^2/4}, \nu=a_5+{a_1a_2/2},
\gamma={a_2^2/4}.\end{equation}

\noindent \textsl{The last two transformations} Let  $X,Y$ and $T$
be the new variables given by
\begin{equation}\label{fortrans}X={x/(4\lambda)},\,\,\,Y={y/(4\lambda^2)},\,\,\,T=t\,
\end{equation} In these new variables the curve writes as
$$TY^2=4X^3+{\mu\over\lambda^2}X^2T+{\nu\over4\lambda^3}XT^2+{\gamma\over4^2\lambda^4}T^3.$$
Finally, with the last change
\begin{equation}\label{fifthtrans}
x={1\over12}{12X\lambda^2+\mu
t\over\lambda^2},\,\,\,y=Y,\,\,\,t=T,\end{equation} the curve is
written in the  Weierstrass' normal form
$$y^2t=4x^3-g_2xt^2-g_3t^3,$$
where
\begin{equation}\label{paramg}{g_2={\mu^2\over12\lambda^4}-{\nu\over4\lambda^3},
g_3=-{\mu^3\over216\lambda^6}+{\mu\nu\over48\lambda^5}-{\gamma\over16\lambda^4}.}\end{equation}

\subsubsection{Proof of Proposition
\ref{propo1}}

\noindent{\textsl{Proof of Proposition \ref{propo1}.}} (i) The
expression of the equation (\ref{Wei}) comes from the
transformations leading ${\cal D}_L$ to its Weierstrass form given
in Section \ref{weiers}. In particular the expressions of the
coefficients of $g_2$ and $g_3$ in terms of $\alpha$, $\beta$ and
$L$ come from formulae (\ref{parama}), (\ref{parammu}) and
(\ref{paramg}).

 (ii) From Proposition \ref{iselliptic}, the curves ${\cal
D}_L$, with $L>L_c$ are elliptic. This means that for $L>L_c$ the
curves $\mathcal{E}_L$ will be also elliptic.

Statement (iii) is a direct consequence of Proposition \ref{l9}. But
it is still necessary to prove that the zero point  $V$ is
transformed into $\widehat{V}=[0:1:0]$. To do this  we should obtain
the image of $V$ under the same transformations as above, but the
image of $V$ under the first transformation, say $V_1:=[1:0:0]$, has
no image for the second transformation. So we use a continuity
argument: Set $V_1(x) = [x:y(x):1]$, where  $y(x)$ is choose in sort
that $V_1(x)$ is on the cubic (\ref{firsteq}), that is, $y(x)$ can
be either $y_+(x)$ or $y_{-}(x)$, where
$$y_{\pm}(x)={\frac {-a_{{1}}x-a_{{2}}\pm\sqrt
{4\,a_{{3}}{x}^{3}+4\,(a_{{4}}+{a_{{1}}}^{2}){x}^{2}+2\,(a_{{1}}
a_{{2}}+4\,a_{{5}})\,x+{a_{{2}}}^{2}} }{2x}}$$

The leading term of $y_{\pm}(x)$ when $x\to +\infty$ is
$\pm\sqrt{a_3\,x}$ (observe that these two numbers are imaginary,
because $a_3 = -1/L^2$). So we have two complex points of the cubic
$[x:y_{\pm}(x):1]$ where $y_{\pm}=\pm\sqrt{a_3\,x}+O(1)$:
$$
V_1(x) = [x:\pm\sqrt{a_3\,x}+O(1):1] = [1:\pm\sqrt{a_3/x}+O(1/x):
1/x]$$ which tends to $V_1=[1:0:0]$ when $x$ tends to $+\infty$.

But if we take the images of these two complex points, in their
first form, by the  transformation (\ref{sectrans}), we obtain
$$V_2(x) = [1/\sqrt{x}:\pm\sqrt{a_3}+O(1/\sqrt{x}):1/\sqrt{x^3}],$$ which tends to
$[0:\pm\sqrt{a_3}:0]=[0:1:0]$. So this point is $V_2$, the image of
$V_1$ obtained by continuity. Now the successive transformations can
be done without any problem and we obtain $\widehat{V}=[0:1:0]$.

 The
point $\widehat{H}$ is the image of $H=[1:0:0]$ under the successive
transformations (\ref{firsttrans}), (\ref{sectrans}),
(\ref{thirdtrans}), (\ref{fortrans}) and (\ref{fifthtrans}).\qed

\subsubsection{Proof of the ellipticity of the curves ${\cal D}_L$ for
$L>L_c$}\label{elmodif}

In this section we prove Proposition \ref{iselliptic}, but first we
notice the following facts:

\smallskip

\noindent \textsl{Fact 1. The curves ${\cal D}_L$ have no singular
points at the infinity of $\C P^2$.} Indeed, setting $ D:=(X+\alpha
T)(Y+\beta T)(T+X+Y)-LXYT,$ we have that
$$\frac{\partial D}{\partial X}(X,Y,0)=Y(Y+2X),\,\frac{\partial D}{\partial
Y}(X,Y,0)=X(2Y+X),$$ so an straightforward computation shows that on
the line $T=0$ the only singular point must satisfy $X=0$ and $Y=0$,
thus not in $\C P^2$.

\smallskip

\noindent \textsl{Fact 2.} An straightforward computation shows that
 the set of affine singular points of ${\cal D}_L$ coincides
with the set of singular points of the invariant function $W$ given
in (\ref{W}). On the other hand, a computation shows that
\begin{equation}\label{p1p2}
  \begin{array}{ll}
    {W}_x' &=(Y+\beta)(X^2-\alpha(1+Y))/(X^2Y)\mbox{ and } \\
    {W}_y' &=(X+\alpha)(Y^2-\beta(1+X))/(XY^2).
  \end{array}
\end{equation}
 \textsl{So for $L>L_c$, the other possible singular point of a curve
${\cal D}_L$ is one of the point intersection of the two parabolas}
$    {\cal P}_1:=\{ X^2=\alpha(1+Y)\}$ \textsl{and} ${\cal P}_2:=
\{Y^2=\beta(1+X)\}.$

\smallskip

\noindent \textsl{Fact 3.} \textsl{The equilibrium point $P=(p,q)$
which is in the positive quadrant, is one of the (at most) four real
points on ${\cal P}_1\cap{\cal P}_2$, and it is on ${\cal D}_L$ and
singular in this curve only for $L=L_c$ (in this case, it is a real
isolated point).}

\smallskip

Now we introduce the some new parameters  $(p,q)$ which correspond
to the coordinates of the fixed point in ${\cal Q}^+$, thus given by
\begin{equation}\label{paramparam}\left\{
  \begin{array}{l}
    p^2=\alpha(1+q), \\
    q^2=\beta(1+p).
  \end{array}
\right.\end{equation} With these new parameters it is easy to obtain
an explicit form for $L_c$:
\begin{equation}\label{Lcexplicit} L_c=\frac{(1+p+q)^3}{(1+p)(1+q)}.
\end{equation}

In fact this explicit expression allows to obtain an straightforward
proof that $L_c>1$ since
$$L_c=(p+q+1)\left(1+\frac{q}{p+1}\right)\left(1+\frac{p}{q+1}\right)>1.$$ Furthermore, after some
computations (skipped here) it is possible to prove that
$$\lim\limits_{(\alpha,\beta)\to(0,0)}(p,q)=(0,0).$$ This implies that
the above one is the sharpest lower bound for $L_c$.

\smallskip

 \noindent \textsl{Fact 4.}  \textsl{Any singular point of ${\cal D}_L$
with $L>L_c$ must be real.} If a singular point is not real, the
conjugate point is also singular and different, and so the (real)
line which support these two points is in the curve, which splits in
a right line and a conic. It is easy to see that, after the
identification of the coefficients the equation of the cubic must be
written as
    $$D=(X+Y+(1-L)T)\left(XY+\beta XT+\alpha YT+\frac{\alpha\beta}{1-L}T^2\right)=0.$$
The point $[0:-1:1]$ is on the cubic, hence $
-L\left(-\alpha+\alpha\beta/(1-L)\right)=0$, which is impossible
because the above number  is positive, since  we have just seen that
$L>L_c>1$.

Taking into account all the above information, to prove Proposition
\ref{iselliptic}, we only need to show that if ${\cal D}_L$ is
singular then $L\leq L_c$.

\begin{lem}\label{nonsingular}
If the curve ${\cal D}_L$ is singular, then $L\leq L_c$.
\end{lem}

\proof Using the Weierstrass form given by (\ref{Wei}) and the
expressions of $g_2$ and $g_3$ given by from formulae
(\ref{parama}), (\ref{parammu}) and (\ref{paramg}). We obtain that
the discriminant of ${\cal D}_L$ is given by
$$
\Delta:=g_2^3-27g_3^2=\frac{1}{4096}\alpha^2\beta^2L^{15}\Delta_1(L;\alpha,\beta),
$$
where $\Delta_1(L)=L^4+\sum_{i=0}^3\delta_{i}(\alpha,\beta)\,L^i$,
and being
$$
\begin{array}{rl}
  \delta_3(\alpha,\beta)= &\alpha \beta-4 \alpha-4 \beta-4, \\
  \delta_2(\alpha,\beta)= &-3 {\alpha}^{2}\beta-3 \alpha {\beta}^{2}+6 {\alpha}^{2}-21
\alpha \beta+6 {\beta}^{2}+4 \alpha+4 \beta+6,
  \\
  \delta_1(\alpha,\beta)=&3 {\alpha}^{3}\beta-21 {\alpha}^{2}{\beta}^{2}+3 \alpha {\beta}^{3
}-4 {\alpha}^{3}+18 {\alpha}^{2}\beta+18 \alpha {\beta}^{2}-4 {
\beta}^{3}+4 {\alpha}^{2}\\
&-25 \alpha \beta+4 {\beta}^{2}+4 \alpha+ 4 \beta-4,
  \\
  \delta_0(\alpha,\beta)=&- \left( \alpha-1 \right) \left( \beta-1 \right)    \left( \alpha+
\beta -1\right) ^{3}.
\end{array}
$$
Obviously $L=0$ gives a singular curve (a union of three lines). To
see that all the real roots of $\Delta_1$ are lower than $L_c$, we
take into account the explicit expression of $L_c$ given by
(\ref{Lcexplicit}), we
 use the change
$$
\left\{    \alpha=\displaystyle{\frac{p^2}{1+q}},\quad
    \beta=\displaystyle{\frac{q^2}{1+p}}\right.,
$$
and also the translation given by $L=\widetilde{L}+L_c$, obtaining
that $\Delta_1(L;\alpha,\beta)=0$ if and only if the numerator of
$\Delta_1(\widetilde{L}+L_c;p^2/(1+q),q^2/(1+p))$ vanishes, which
happens if and only if $
\widetilde{\Delta}_1(\widetilde{L}):=\widetilde{L}\,\left[\sum_{i=0}^3
d_i(p,q) \widetilde{L}^i \right]=0$, where
$$
\begin{array}{rl}
  d_3(p,q)=&(1+p)^3(1+q)^3, \\
  d_2(p,q)=&(1+p)^2(1+q)^2\left({p}^{2}{q}^{2}+12 {p}^{2}q+12 p{q}^{2}+8 {p}^{2}+20 pq+8 {q}^{2}+
8 p+8 q
\right),\\
  d_1(p,q)=&(1+p)^2(1+q)^2\left(9 {p}^{4}{q}^{3}+9 {p}^{3}{q}^{4}+60 {p}^{4}{q}^{2}+93 {p}^{3}{q}^
{3}+60 {p}^{2}{q}^{4}+64 {p}^{4}q+228 {p}^{3}{q}^{2}\right.\\
{}&\left.+228 {p}^{2}{q }^{3}+64 p{q}^{4}+16 {p}^{4}+176
{p}^{3}q+312 {p}^{2}{q}^{2}+176 p{q}^{3}+16 {q}^{4}+32 {p}^{3}+160
{p}^{2}q\right.\\
&\left.+160 p{q}^{2}+32 {q}^{
3}+16 {p}^{2}+48 pq+16 {q}^{2}\right),\\
 d_0(p,q)=&pq \left( {p}^{2}+pq+{q}^{2}+p+q \right)  \left( 3 pq+4 p+4 q+4
 \right) ^{3}.
\end{array}
$$
Since $\widetilde{\Delta}_1$ is a polynomial with positive
coefficients, its real roots are negative or zero, and thus,  the
real roots of $\Delta_1$ are less or equal than $L_c$.\qed

\section{Behavior of $\theta_{b,a}(h)$.}\label{behav}

To prove Theorem \ref{periodsinI}, as well as to study  the periods
appearing in the family $F_{b,a}$ we need to analyze the rotation
number function $\theta_{b,a}(h)$, with special emphasis on  its
asymptotic behavior  when $h\to h_c^+$,  and when $h\to+\infty$.
Regarding this last question, in \cite{CGM11a}, there were pointed
out strong numerical evidences that $\lim\limits_{h\to +\infty}
\theta_{b,a}(h)=2/5$. Notice that these evidences agree with the
fact, proved in \cite{BR}, that for the autonomous case the rotation
number associated to the invariants curves (\ref{elipautonom}), say
$\theta_a(h)$, satisfies $\lim\limits_{h\to +\infty}
\theta_{a}(h)=1/5$ so it is  independent of the value of the
parameter $a$. In this section we prove

\begin{propo}\label{p1} The following statements hold
\begin{itemize}
  \item[(i)] The rotation number map $\theta_{b,a}(h)$ is analytic  on
  $(h_c,+\infty)$.
      \item[(ii)] $\lim\limits_{h\to +\infty}
\theta_{b,a}(h)=\displaystyle{\frac{2}{5}}$.
  \item[(iii)] The rotation number map $\theta_{b,a}(h)$ is continuous in   $[h_c,+\infty)$. Furthermore, setting
  \begin{equation}\label{funciosigma}
  \sigma(a,b):=\lim\limits_{h\to h_c^+}\theta_{b,a}(h),
   \end{equation}
\noindent it holds that \begin{eqnarray} \sigma(a,b)
&=&\frac{1}{2\pi}\arccos\left(\frac{1}{2}\left[-1-\displaystyle{\frac{a+b\,x_c}{x_c\,(b+x_c)}}\right]\right)
\label{sigma1}\\
&=&\displaystyle{ \frac{1}{2\pi}}\arccos\left(\displaystyle{
\frac{1}{2}}\left[-2+\displaystyle{\frac{1}{x_c\,y_c}}\right]\right).
  \label{sigma2}
  \end{eqnarray}

\end{itemize}
\end{propo}

To prove the statements (i) and (ii) of Proposition \ref{p1}, it
will be helpful to use the conjugation of $F_{b,a}$ given by the map
$\widehat{G}$ defined in (\ref{accio}), for which its invariant
level sets are described by the Weierstrass normal form ${\cal E}_L$
defined in (\ref{Wei}), and then to apply a similar machinery as in
\cite{BR}. First, we  recall that given the elliptic curve ${\cal
E}_L$, for $L>L_c$, there exist two positive numbers $\omega_1$ and
$\omega_2$ depending on $\alpha,\beta$ and $L$ and a lattice in
$\mathbb{C}$
$$
\Lambda=\{2n\omega_1+2m\,i\omega_2\,\mbox{ such that } (n,m)\in
\mathbb{Z}^2\}\subset\mathbb{C},
$$
such that the \emph{Weierstrass $\wp$ function relative to}
$\Lambda$
$$
\wp(z)=\frac{1}{z^2}+\sum\limits_{\lambda\in\Lambda\setminus\{0\}}\left[\frac{1}{(z-\lambda)^2}-\frac{1}{\lambda^2}\right]
$$
gives a parametrization of ${\cal E}_L$. This is because the map
\begin{equation}\label{phi}
    \begin{array}{cccl}
      \phi: & \mathbb{C}/\Lambda & \longrightarrow & {\cal E}_L\\
      {} & z & \longrightarrow & \left\{
                                   \begin{array}{cl}
                                     [\wp(z):\wp'(z):1] & \mbox{if}\quad z\notin\Lambda, \\
                                    \phantom{ }[0:1:0]=\widehat{V} & \mbox{if}\quad z\in\Lambda,
                                   \end{array}
                                 \right.
    \end{array}
\end{equation}
is an holomorphic homeomorphism, and therefore
\begin{equation}\label{diffeq}
\wp'(z)^2=4\wp(z)^3-g_2\wp(z)-g_3
\end{equation}
(see \cite[Lemma 5.17]{K} and \cite[Proposition 3.6]{S}). Another
interesting property of the parametrization of ${\cal E}_L$ given by
$\phi$ is that it is injective for the real interval $(0,\omega_1)$
onto the real unbounded branch of ${\cal E}_L$ whose points have
\emph{negative $y$--coordinates}, and that
$\phi(\omega_1)=[e_1,0,1]$ or, in other words:
 \begin{equation}\label{e1wp}
e_1=\wp(\omega_1)
 \end{equation}
(see again \cite{K,S}). Finally, observe that
$\phi(0)=\phi(2\omega_1)=\widehat{V}$, and $\lim\limits_{u\to
0}\wp(u)=+\infty$. Taking all the preceding considerations into
account, and by direct integration of the differential equation
(\ref{diffeq}) on $[0,u)$, we have that in real variables
 \begin{equation}\label{up}
u=\displaystyle{\displaystyle{\int}_{\wp(u)}^{+\infty}
\frac{\mathrm{d}s}{\sqrt{4s^3-g_2s-g_3}}}.
 \end{equation}

\noindent{{\sl Proof of Proposition \ref{p1}  (i).}  Observe that
equation (\ref{Wei}) can be written as
$$ty^2=4(x-te_1)(x-te_2)(x-te_3),$$ where
$$e_1+e_2+e_3=0,\,\,\,e_1e_2+e_2e_3+e_3e_1=-g_2/4,\,\,\,e_1e_2e_3=g_3/4,$$
and $e_3<e_2<e_1$.  Then, equation (\ref{up}) writes:
 \begin{equation}\label{up2}
u=\displaystyle{\displaystyle{\int}_{\wp(u)}^{+\infty}
\frac{\mathrm{d}s}{\sqrt{4(s-e_1)(s-e_2)(s-e_3)}}},
 \end{equation}
and from the above relation and equation (\ref{e1wp})  we get:
 $$
\omega_1=\displaystyle{\int}_{e_1}^{+\infty}{{\rm
d}s\over\sqrt{4(s-e_1)(s-e_2)(s-e_3)}},
 $$

From Proposition \ref{propo1}, the action of $\widehat{G}$ on each
level ${\cal E}_L$, with $L>L_c$ is conjugate to a rotation of angle
$2\pi\Theta(L)$, where $\Theta$ is the rotation number associated to
$\widehat{G}_{|{\cal E}_L}$. On the other hand from the expression
of $\widehat{H}$ given in (\ref{eqhtilde}) we know that its
$y$--coordinate is negative, so from the above considerations on the
parametrization $\phi$ we know that there exists $u\in(0,\omega_1)$
such that it has the following representation
$$\exp\left(\frac{2\pi i}{2\omega_1}u\right)=\exp\left(2\pi i\Theta(L)\right)\in[0,2\omega_1]/\Lambda\cong \mathbb{S}^1,$$
so $u=2\omega_1\Theta(L)$, and since $\widehat{H}$ is the image of
$\widehat{V}$ under the transformation $\widehat{G}_{|{\cal E}_L}$,
it is given by $\phi(2\omega_1 \Theta(L))$, and therefore
$$\wp(2\omega_1 \Theta(L))=X(L),$$ where $X(L)$ is the abscissa of
point $\widehat{H}$ given in (\ref{eqhtilde}). Applying (\ref{up2})
 we obtain that
$$
2\omega_1\Theta(L)=\displaystyle{\int}_{X(L)}^{+\infty}{{\rm
d}s\over\sqrt{4(s-e_1)(s-e_2)(s-e_3)}},
$$
hence
\begin{equation}\label{2theta}2\Theta(L)=\displaystyle{\displaystyle{\int}_{X(L)}^{+\infty}{{\rm
d}s\over\sqrt{(s-e_1)(s-e_2)(s-e_3)}}\over
\displaystyle{\int}_{e_1}^{+\infty}{{\rm
d}s\over\sqrt{(s-e_1)(s-e_2)(s-e_3)}}}.\end{equation}

By using the change of variables $s=e_1+1/r^2$ and
$r\sqrt{e_1-e_3}=u$, it can be proved that equation (\ref{2theta})
becomes
\begin{equation}\label{eq17}
2\Theta(L)=\displaystyle{\displaystyle{\int}_{0}^{\sqrt{\frac{e_1-e_3}{\nu}}}{{\rm
d}u\over\sqrt{(1+u^2)(1+\varepsilon u^2)}}\over
\displaystyle{\int}_{0}^{+\infty}{{\rm
d}u\over\sqrt{(1+u^2)(1+\varepsilon u^2)}}},
\end{equation}
where $\nu=X(L)-e_1$, and  $\varepsilon=(e_1-e_2)/(e_1-e_3)$.

The analyticity of $\Theta(L)$ in $(L_c,+\infty)$ can be derived now
from expression (\ref{eq17}),  by taking into account that when
$L>L_c$, $e_1$, $e_2$, $e_3$ are different, and that from the
geometry of the problem we have $e_3<e_2<e_1<X(L)$, so all the
elements in the left hand side of (\ref{eq17}) are holomorphic in a
neighborhood of $(L_c,+\infty)$. Hence $\Theta(L)$ is analytic on
$(L_c,+\infty)$.\qed

\smallskip

Observe that $\lim\limits_{h\to +\infty}
\theta_{b,a}(h)=\lim\limits_{L\to +\infty} \Theta(L)$, hence our
main objective now is to prove that $\lim\limits_{L\to +\infty}
\Theta(L)=2/5$. To do this  we will need the following auxiliary
result

\begin{lem}\label{lem4BR}(\cite[Lemma 4]{BR})
Let $\lambda,\varepsilon,\gamma$ be positive numbers. For any map
$\phi(\varepsilon)$ such that $\lim\limits_{\varepsilon\to
0}\phi(\varepsilon)=0,$  and $\lambda+\phi(\varepsilon)>0$, set
$$
N(\varepsilon,\lambda,\gamma)=\displaystyle{\displaystyle{\int}_{0}^{\frac{\lambda+\phi(\varepsilon)}{\varepsilon^\gamma}}}{{\rm
d}u\over\sqrt{(1+u^2)(1+\varepsilon u^2)}}, \mbox{ and }
D(\varepsilon)=\displaystyle{\displaystyle{\int}_{0}^{+\infty}{{\rm
d}u\over\sqrt{(1+u^2)(1+\varepsilon u^2)}}}.
$$
Then $D(\varepsilon)\sim (1/2)\ln(1/\varepsilon)$, and if
$\gamma<1/2$ we have $N(\varepsilon,\lambda,\gamma)\sim
\gamma\ln(1/\varepsilon)$, where $\sim$ denotes the equivalence with
the leading term of the  asymptotic development at zero.
\end{lem}

\noindent{{\sl Proof of Proposition \ref{p1} (ii).} Following the
steps in \cite{BR}, we will compute
$\lim\limits_{L\to+\infty}\Theta(L)$ by using the asymptotic
developments of the elements involved in equation (\ref{eq17}).

Solving the cubic equation $4\,s^3-g_2\,s-g_3=0$ by leaving the
solutions in terms of $g_2$ and $g_3$;  using  the expressions of
these coefficients in terms of $\alpha$, $\beta$ and $L$ given by
Proposition \ref{propo1}; and after some computations done with a
symbolic algebra system (Maple 12), we obtain that
$$
e_1\sim L^4/48,\, e_2\sim L^4/48\, \mbox{ and } e_3\sim -L^4/24,
$$
when $L$ tends to infinity, and
$$
e_1-e_2\sim \alpha \beta\, L^{3/2}\mbox{ and } e_1-e_3\sim L^4/16.
$$
Setting $\varepsilon=(e_1-e_2)/(e_1-e_3)$, we obtain
$\varepsilon\sim 16\,\alpha\beta/L^{5/2}$ (when $L$ tends to
infinity), so

\begin{equation}\label{Lepsilon}
L\sim\displaystyle{\frac{2\,2^{3/5}\,(\alpha\beta)^{2/5}}{\varepsilon^{2/5}}}.
\end{equation}
(when $\varepsilon$ tends to zero). From the expression of
$\widehat{H}$ given in (\ref{eqhtilde}) we have
$$
\nu=X(L)-e_1\sim \alpha\beta\,L^2/4,
$$
and therefore, using (\ref{Lepsilon}), we get
$$
\displaystyle{\sqrt{\frac{e_1-e_3}{\nu}}}\sim\displaystyle{\frac{L}{2\sqrt{\alpha\beta}}}\sim
A\,\varepsilon^{-2/5},$$where $A=2^{3/5}\,(\alpha\beta)^{-1/(10)}$.
So from equation (\ref{eq17}), and Lemma \ref{lem4BR} we finally
obtain
$$
2\Theta(L)=\displaystyle{\frac{N(\varepsilon,A,2/5)}{D(\varepsilon)}}\sim
\displaystyle{\frac{\frac{2}{5}
\ln(1/\varepsilon)}{\frac{1}{2}\ln(1/\varepsilon)}}=\frac{4}{5},
$$
and therefore $\lim\limits_{L\to\infty}\Theta(L)=2/5$.\qed

It is well known (see \cite{BR2} for instance, and also
\cite{CGM07,Z}) that under some regularity conditions the rotation
number function of a differentiable planar map $F$ near an elliptic
fixed point, can be extended continuously to the fixed point. In
this sense, the proof of statement (ii) follows by applying the
following version of \cite[Proposition 8]{BR2}.

\begin{lem}\label{Zdm} Let $U$ be an open set of $\R^2$ and a $\mathcal{C}^1$ map
$F:U\rightarrow U$. Suppose that $F$ has a unique elliptic fixed
point $P$ in $U$ (that is $DF(P)$ has two distinct non real
eigenvalues with modulus $1$). Suppose that $F$ has a first integral
$V$ in $U$, attaining a strict minimum $h_p:=V(P)$ at the point $P$,
and such that
\begin{itemize}
  \item[(i)] For $h>h_p$, the invariant level sets $\{V=h\}$ (except the level set
$\{V(P)=h_P\}$) are closed curves, surrounding $P$.
\item[(ii)] For $h>h_p$,  the sets $\overline{\mathrm{Int}\left(\{V=h\}\right)}$
are starlike with respect to $P$, and each half line from $P$ cuts
$\{V=h\}\cap U$ exactly in one point.
  \item[(iii)] On each curve $\{V=h\}$ for $h>h_p$,  the map  $F_{|\{V=h\}}$ is conjugate to a rotation
of angle $2\pi\theta(h)$.
\end{itemize}
Then,
\begin{equation}\label{expZdm}
\theta(P):=\lim\limits_{h\to h_P}
\theta(h)=\displaystyle{\frac{1}{2\pi}}
\arccos\left(\frac{\mathrm{Trace}(DF(P))}{2}\right).
\end{equation}
\end{lem}

\noindent{\sl Proof of Proposition \ref{p1}  (iii).} As a direct
consequence of Lemma \ref{Zdm}, the continuous extension of
$\theta_{b,a}(h)$ to $h_c^+$, is
$$
\sigma(a,b)=\displaystyle{\frac{1}{2\pi}\arccos\left(\displaystyle{\frac{1}{2}}\mathrm{Trace}(DF_{b,a}(x_c,y_c))\right)}.
$$
 Then, the  expressions
(\ref{sigma1}--\ref{sigma2}) are then easily derived from
(\ref{puntfix}) after some easy manipulations.\qed

\section{Periodic solutions}\label{provamainresults}

\subsection{The rotation interval $I(a,b)$. Proof of Theorem \ref{periodsinI}}\label{provateo1}

In this Section we prove Theorem \ref{periodsinI}, and we give a
more precise characterization of the rotation number interval
$I(a,b)$.

\noindent \textsl{Proof of Theorem \ref{periodsinI}.} In proposition
\ref{p1}, it is proved that $\theta_{b,a}(h)$ is an analytic
function  on $(h_c,+\infty)$ and continuous in $[h_c,+\infty]$; that
$\lim\limits_{h\to h_c^+} \theta_{b,a}(h)=\sigma(a,b)$; and
$\lim\limits_{h\to+\infty}\theta_{b,a}(h)=2/5$.

It is clear then that for any $\theta$ in the interval $I(a,b)$
defined by
$$I(a,b):=\left\langle \lim\limits_{h\to h_c^+}
\theta_{b,a}(h),\lim\limits_{h\to +\infty} \theta_{b,a}(h)
\right\rangle=\left\langle \sigma(a,b),{2\over5} \right\rangle,$$
where $\langle c,d\rangle=(\min(c,d),\max(c,d))$, there is a value
of $h\in[h_c,+\infty)$ such that the action of $F_{b,a}$ restricted
to $\mathcal{C}_h^+$ is a rotation with rotation number $\theta$,
since $$I(a,b)\subseteq
\mathrm{Image}\left(\theta_{b,a}((h_c,+\infty))\right).$$\qed

To have a better description of the interval $I(a,b)$, it would be
interesting to know the relative position of $\sigma(a,b)$ with
respect to $2/5$ for given values of $a$ and $b$. This is issue is
considered in Proposition \ref{propocondiciosigma} and Corollary
\ref{corolgamma} (see also Figure 1), but before proving them we
need a preliminary result.

Set $J(a,b):=
\mathrm{Image}\left(\theta_{b,a}\left((h_c,+\infty)\right)\right)\subset\big[0,{1\over2}\big].$
As as stated above, it holds  that  $I(a,b)\subseteq J(a,b)$.
Observe that both intervals are the same when the function
$\theta_{b,a}(h)$ is monotonic, but notice that in \cite{CGM11a} it
is proved that there are values of $a$ and $b$ for which the map is
non--monotonic. With respect the interval $J(a,b)$, we prove:

\begin{lem}\label{lema123}
The following statements hold.
\begin{itemize}
  \item[(i)] For any $a,b>0$,  $\sigma(a,b)>\displaystyle{\frac{1}{3}}$.
  \item[(ii)]
  $\lim\limits_{b\to+\infty}\sigma(b^2,b)=\displaystyle{\frac{1}{2}}$. Therefore, when $b$
  is sufficiently large  $\sigma(b^2,b)>2/5$, and then $I(b^2,b)=\left(2/5,\sigma(b^2,b)\right)$.
  \item[(iii)]  $\lim\limits_{b\to
  0^+}\sigma(b^2,b)=\displaystyle{\frac{1}{3}}$. Therefore, when $b$
  is small enough $\sigma(b^2,b)<2/5$, and then $I(b^2,b)=\left(\sigma(b^2,b),2/5\right)$.
\end{itemize}
\end{lem}

As a first direct consequence of the above result, we have

\begin{propo}\label{unionofJ} It holds that
$$\left({1\over3},{1\over2}\right)\subset \bigcup_{a>0,\,b>0}J(a,b).$$ Furthermore
$$\left({1\over3},{1\over2}\right)\subset \bigcup_{b>0}J(b^2,b).$$
\end{propo}

\noindent \textsl{Proof of Lemma \ref{lema123}.} (i) Using that
$$
\sigma(a,b)=\displaystyle{\frac{1}{2\pi}\arccos\left(\displaystyle{\frac{1}{2}}\mathrm{Trace}(DF_{b,a}(x_c,y_c))\right)},
$$
we have that $\sigma(a,b)>1/3$ if and only if
$\mathrm{Trace}(DF_{b,a}(x_c,y_c))/2<\cos(2\pi/3)=-1/2$ or,
equivalently using equation (\ref{sigma1}),
$$
-1-\displaystyle{\frac{a+b\,x_c}{x_c\,(b+x_c)}}<-1,
$$
which is true, since $a$, $b$ and $x_c$ are positive.

\noindent (ii) From equation (\ref{puntfix}) we have that
$y_c>\sqrt{b}$. Suppose that $a=b^2$, then
$x_c=\sqrt{a+y_c}>\sqrt{a}=b.$ So $\lim\limits_{b\to +\infty}
x_c=\lim\limits_{b\to +\infty} y_c=+\infty$ and then, using equation
(\ref{sigma2})
$$\lim\limits_{b\to+\infty}\theta(b^2,b)=\lim\limits_{b\to+\infty}\displaystyle{ \frac{1}{2\pi}}\arccos\left(\displaystyle{
\frac{1}{2}}\left[-2+\displaystyle{\frac{1}{x_c\,y_c}}\right]\right)=\displaystyle{
\frac{1}{2\pi}}\arccos(-1)=\displaystyle{\frac{1}{2}}.$$

\noindent (iii) When $a=b^2$, we consider the polynomials
$R(x,y):=x^2-y-b^2$ and $S(x,y):=y^2-x-b$ whose common zeros are the
fixed points of $F_{b,a}$. Observe that $y_c$ must be a root of the
resultant
$$
\mathrm{Resultant}(R,S;x)=-y(-y^3+2by+1),
$$
so $y_c^3=2by_c+1$, thus $y_c>1$ and $y_c<y_c^3=2by_c+1$, hence
$$
y_c<\frac{1}{1-2b}<2\,\mbox{ if } b<\frac{1}{4}.
$$
So if $b<1/4$ we have
$$
1<y_c^3=1+2by_c<1+4b;
$$
since the right hand side of the above inequality tends to $1$ when
$b$ tends to $0$, we obtain that $y_c$ tends to $1$. Since
$x_c=y_c^2+b$, we obtain that $\lim\limits_{b\to 0^+}
x_c=\lim\limits_{b\to 0^+} y_c=1$, hence
$$\lim\limits_{b\to 0^+}\theta(b^2,b)=\lim\limits_{b\to 0^+}\displaystyle{ \frac{1}{2\pi}}\arccos\left(\displaystyle{
\frac{1}{2}}\left[-2+\displaystyle{\frac{1}{x_c\,y_c}}\right]\right)=\displaystyle{
\frac{1}{2\pi}}\arccos(-1)=\displaystyle{\frac{1}{2}}.$$ \qed

The next result characterizes the relative position of $\sigma(a,b)$
with respect $2/5$.

 \begin{propo}\label{propocondiciosigma} Let $(x_c,y_c)$
be the fixed point of the map $F_{b,a}$ given by (\ref{puntfix}).
Set $R(x,y):=\phi\, y\,-\,b\,x-a$. Then,  $\sigma(a,b)<2/5$ if and
only if $R(x_c,y_c)>0$; $\sigma(a,b)=2/5$ if and only if
$R(x_c,y_c)=0$; and $\sigma(a,b)>2/5$ if and only if $R(x_c,y_c)<0$.
\end{propo}

\proof Observe that using expression (\ref{sigma1}), the condition
$\sigma(a,b)\leq 2/5$ gives
$$
\frac{1}{2}\left[-1-\displaystyle{\frac{a+b\,x_c}{x_c\,(b+x_c)}}\right]\geq
\cos\displaystyle{\left(\frac{4\pi}{5}\right)}=-\displaystyle{\frac{1+\sqrt{5}}{4}},
$$
thus
$$
\displaystyle{\frac{a+b\,x_c}{x_c\,(b+x_c)}}\leq\phi-1=\displaystyle{\frac{1}{\phi}},
$$
where $\phi=(1+\sqrt{5})/2$. By using (\ref{puntfix}) this last
inequality turns to
$$
\phi\,(a+b\,x_c)\leq x_c\,(x_c+b)=x_c^2+b\,x_c=a+y_c+b\,x_c.
$$
Using again $\phi-1=1/\phi$, we obtain that $ \phi y_c-bx_c-a\geq
0.$\qed

\medskip

In order to obtain a characterization of the condition
$\sigma(a,b)=2/5$ in terms of the parameters $a$ and $b$ instead of
$x_c$ and $y_c$, we isolate $a$ and $b$ from the equations
(\ref{puntfix}) and plug them in the equation  $R(x_c,y_c)=0$,
obtaining
$$
y_c(\phi -x_cy_c+1)=0.
$$
Using that $\phi+1=\phi^2$, we obtain that $R(x_c,y_c)=0$ if and
only if $y_c=\phi^2/x_c$, thus taking the parameter $t:=x_c$ we have
that
$$
  \displaystyle{a=\frac{t^3-\phi^2}{t}}\, \mbox{ and }
  \displaystyle{b=\frac{\phi^4-t^3}{t^2}}.$$

A simple computation shows that the above formulae describe a curve
with a single branch in the parameter's space
$\mathcal{P}:=\{(a,b),\,a,b>0\}$ when
$\phi^{\frac{2}{3}}<t<\phi^{\frac{4}{3}}$. In summary, we have
proved the following corollary:
\begin{corol}\label{corolgamma}
The curve $\sigma(a,b)=2/5$ for $a,b>0$ is given by
$$\Gamma:=\left\{\sigma(a,b)=2/5,\,a,b>0\right\}=
\left\{(a,b)=\left(\frac{t^3-\phi^2}{t},\frac{\phi^4-t^3}{t^2}\right),
\,t\in(\phi^{\frac{2}{3}},\phi^{\frac{4}{3}})
\right\}\subset\mathcal{P}.$$
\end{corol}

It is easy to check that this curve cuts the axes at the points
$(0,\sigma_*)$  and $(\sigma_*,0)$ where
$\sigma_*=\phi^{\frac{5}{3}}\simeq 2.2300$, and it splits the
parameter space into two connected components. By using Lemma
\ref{lema123} (ii) and (iii), it is easy to obtain that one of them
corresponds to the set $\{\sigma(a,b)<2/5\}$ and the other to
$\{\sigma(a,b)>2/5\}$, and that they have the geometry depicted in
Figure~1.

\begin{center}
\includegraphics[scale=0.45]{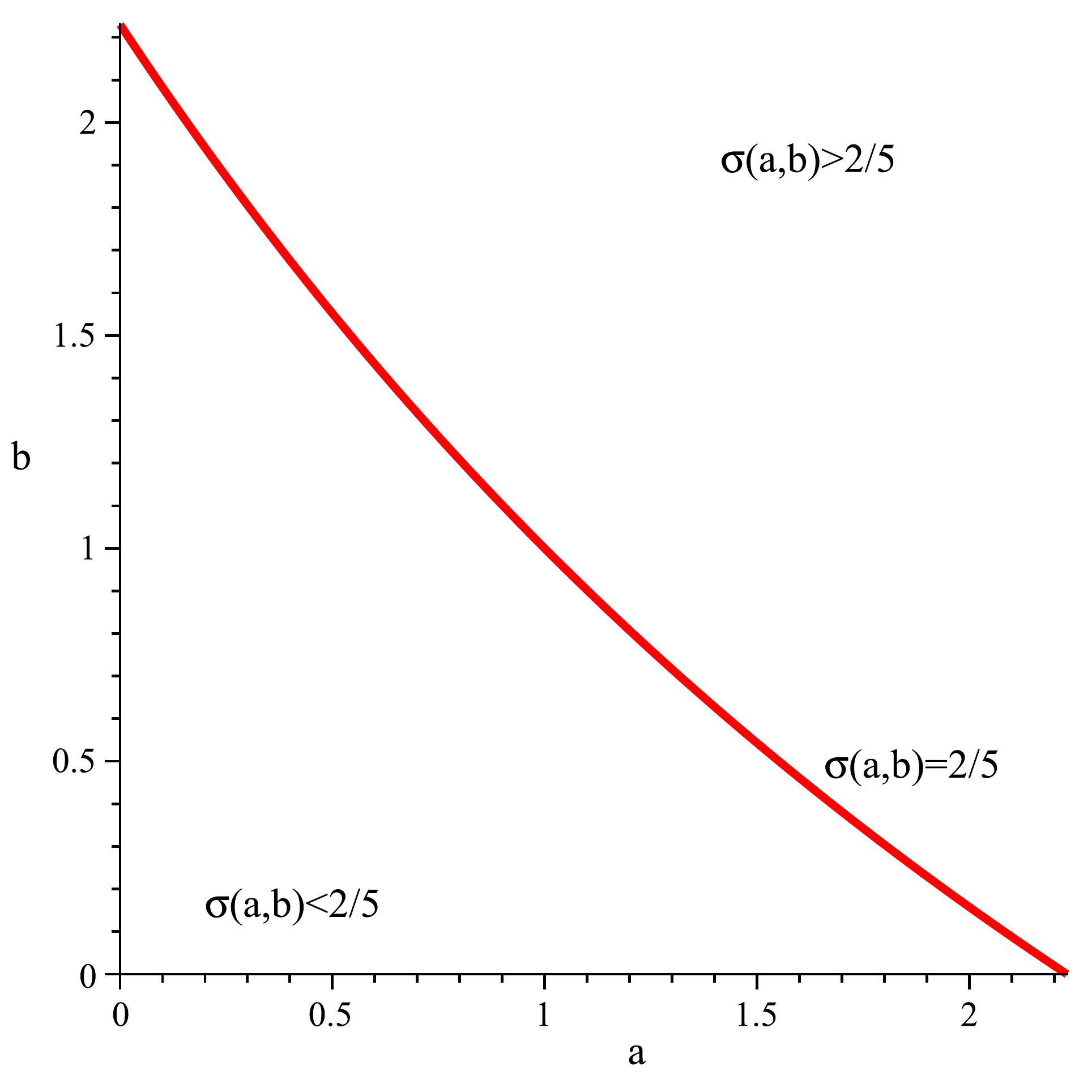}
\end{center}
\begin{center}
Figure 1: The sets $\{\sigma(a,b)<2/5\}$ and $\{\sigma(a,b)>2/5\}$,
and the curve $\sigma(a,b)=2/5$.
\end{center}

\subsection{Forbidden periods and periods not appearing in ${\cal Q}^+$}\label{forbidden}

\begin{propo}\label{propoimp1}
The map $F_{b,a}$ do not have periodic orbits of prime periods $2$
and $3$.
\end{propo}

\proof Observe that the orbits of the dynamical system defined by
the iterations of $F_{b,a}$ are also given by the planar first order
system of difference equations given by $$\left\{
    \begin{array}{ll}
      u_{n+1}u_n=a+v_n, \\
     v_{n+1}v_n=b+u_{n+1}.
    \end{array}\right.$$

If the above system has a periodic solution of prime period $2$,
then we have that $u_{n+1}u_n$ is constant, so it also must be
$a+v_n$, and therefore $v_n$ and $u_n$ are constant.

If there is a periodic solution of prime period $3$ of the above
system we have
\begin{equation}\label{est1}
u_{n+2}u_{n+1}u_n=c=u_{n+2}(a+v_n),
\end{equation}
where $c$ is a constant; and also we have
$$
v_{n+2}v_{n+1}v_n=d=(b+u_{n+2})v_{n},
$$
where $d$ is a constant. Hence
$$
u_{n+2}=\frac{d}{v_n}-b.
$$
Using equation (\ref{est1}) we have
$c=(a+v_n)\left(\frac{d}{v_n}-b\right)$, so
$$
b\,v_n^2\,+\,(ab+c-d)\,v_n\,-\,ad\,=\,0,
$$
and therefore $v_n$ is constant.\qed

\begin{propo}\label{propoimp4}
The map $F_{b,a}$ do not have periodic orbits of prime periods $4$,
$6$ and $10$
 in ${\cal Q}^+$.
\end{propo}

 Observe, however,  that there can exist periodic orbits of minimal
period $4$, $6$ and $10$ in $\R^2\setminus {\cal Q}^+$ for some
values of $a$ and $b$. For instance, if $a=1,b=2,$ then the curves
${\cal C}_h$ with  $h=3$; $h=4+\sqrt{13}$; and $h\simeq
9.8309187775$ gives rise to continua of initial conditions with
minimal period $4$, $6$ and $10$ respectively. These curves are
located in $\mathbb{R}^2\setminus{\cal Q}^+$. Of course, the period
$10$ appears when considering the map $F_{1,1}$, but not as a prime
period.

\smallskip

\noindent \textsl{Proof of Proposition \ref{propoimp4}.} We consider
again the coordinates $X$, $Y$ and $T$  introduced in
(\ref{changeparameters}), and the curves ${\cal D}_L$. Observe that,
there exist periodic orbits with minimal period $2n$ in the curve
$\mathcal{D}_L$ if and only if $2nH=V$, or in other words, if and
only if $nH=-nH$. Using that $V*V=R:=[-\alpha:0:1]$, we have
$$nH=V\Longleftrightarrow nH=-nH=nH*(V*V) \Longleftrightarrow nH=nH*R.$$
The last of the above condition means that there exist a $2n$
periodic orbit if and only if  $R$ belongs to the tangent line to
${\cal D}_L$ in $nH$. In the following we will see that for periods
$4$, $6$ and $10$ this condition implies  $L<L_c$, thus proving that
these periods do not appear in ${\cal Q}^+$.

\smallskip

\noindent \textsl{{Period 4.}} Since $2H=(H*H)*V=[0:-1:1]=:Q$, there
exist a $4-$periodic orbit if and only if  $R$ belongs to the
tangent line to ${\cal D}_h$ in $Q$. An straightforward computation
allows us to find the explicit expression of the tangent line to
$2H$: $Y=m(L;\alpha,\beta)\,X\,+\,n(L;\alpha,\beta)$, where $$
m(L;\alpha,\beta)=
\displaystyle{-\frac{L+\alpha(\beta-1)}{\alpha(\beta-1)} }\,\mbox{
and }\, n(L;\alpha,\beta)= -1.$$ This tangent line passes through
$R$ if and only if
\begin{equation}\label{condicioR}
-\alpha\,m(L;\alpha,\beta)+n(L;\alpha,\beta)=0,
\end{equation}
or equivalently, after some manipulations, if and only if $L+ \left(
\beta-1 \right)  \left(\alpha-1 \right)=0$. That is, when
$$
L=L_4:=-\left( \beta-1 \right)  \left(\alpha-1 \right).
$$

 To see that $L_4<L_c$, we consider again the
parameters $p$ and $q$ given by the coordinates of the fixed point
in ${\cal Q}^+$, given by (\ref{paramparam}), that is,  we  use the
change
\begin{equation}\label{alphatopq}
\left\{    \alpha=\displaystyle{\frac{p^2}{1+q}},\quad
    \beta=\displaystyle{\frac{q^2}{1+p}}\right.,
\end{equation}
obtaining that
$$
L_4={\frac { \left( -{q}^{2}+p+1 \right)  \left( {p}^{2}-q-1 \right)
}{
 \left( 1+p \right)  \left( 1+q\right)
 }}<\frac{(1+p+q)^3}{(1+p)(1+q)}=L_c.
$$
The last inequality can be obtained by subtraction of the
numerators, for instance.

\smallskip

\noindent \textsl{{Period 6.}} A computation shows that
$$
3H=[-L-\alpha\, \left( \beta-1 \right) : L\beta+\beta\, \left(
\beta-1
 \right)  \left( \alpha-1 \right) :\beta-1]
$$

Proceeding as in the period $4$ case, we compute the tangent line to
${\cal D}_L$ at $3H$, and we impose condition (\ref{condicioR}),
obtaining that there exists a continuum of periodic orbits
characterized by the energy level $L$ if and only if either
$L=0<L_c$ or if $L$ is a root of
$$P_6(L):={L}^{2}+ \left(3 \alpha \beta-2 \alpha -2 \beta+1 \right) L+
 \left( \beta-1 \right)  \left( \alpha-1 \right)  \left( 2 \alpha
\beta-\alpha-\beta +1\right).$$ Using the change (\ref{alphatopq})
and the translation $L=\widetilde{L}+L_c$ we obtain that $P_6(L)=0$
if and only if $\widetilde{L}$ is not a root of
$$\begin{array}{rl}
 \widetilde{P}_6(\widetilde{L}):= & \left( 1+q \right) \left( 1+p \right) \widetilde{L}^{2}+  \left(3
{p}^{2}{q}^{2}+6 {p}^{2}q+6 p{q}^{2}+4 {p}^{2}+13 pq+4 {q}^{2} +7
p+7 q+3\right)
 \widetilde{L} \\
     &+  \left( 1+ q
\right) ^{2} \left( 1+p \right) ^{2} \left( 2 pq+3 p+3 q+3
 \right).
  \end{array}$$
Observe that all the coefficients of
$\widetilde{P}_6(\widetilde{L})$ are positive when $p$ and $q$ are
positive, so all its real roots are negative. Therefore, all the
real roots of $P_6$ are less than $L_c$.

\smallskip

\noindent \textsl{{Period 10.}} A computation done with the aid of a
computer algebra system gives $5H:=[A:B:C]$, where
$$\begin{array}{rl}
    A:=&-\left[\left( \alpha\,\beta -1\right) L+\alpha\,\beta\, \left( \beta-1
 \right)  \left(\alpha-1 \right)
\right]\left[{L}^{2}+\alpha \left( \beta-1 \right)  \left( \beta-2
\right) L-
\alpha \left( \beta-1 \right) ^{3} \left( \alpha-1\right)\right], \\
    B:=&-\beta\left[L+ \left( \beta-1 \right)  \left(\alpha-1 \right) \right]  \left[ L+\alpha \left( \beta-1 \right)  \right]
    \left[
{L}^{2}+
 \left( -2\alpha+1-2\beta+3\alpha\beta \right) L+ \right.\\
 &\left.+\left( \beta
-1 \right)  \left( \alpha-1\right)  \left(
2\alpha\beta-\alpha- \beta+1 \right)  \right],\\
C:=&\left[ \beta L+ \left( \beta-1 \right)
 \left( \beta(\alpha-1)+1 \right)  \right]  \left[ L+ \left(
\beta-1 \right)  \left( \alpha-1\right)  \right]  \left[  \left(
\alpha\beta-1 \right) L+\alpha\beta \left( \beta-1 \right)
 \left( \alpha-1\right)  \right].
  \end{array}
$$
We compute the tangent line to ${\cal D}_L$ at $5H$ and we impose
condition (\ref{condicioR}), obtaining that there exists a continuum
of periodic orbits characterized by the energy level $L$ if and only
if
$$ P_{10,a}(L;\alpha,\beta)\, P_{10,b}(L;\alpha,\beta)\,
P_{10,c}(L;\alpha,\beta)=0,$$ where $$
P_{10,a}(L;\alpha,\beta)=\left( \alpha\,\beta-1 \right)
L+\alpha\,\beta\, \left( \beta-1
 \right)  \left( -1+\alpha \right),\, P_{10,b}(L,\alpha,\beta)=\left(L+\alpha(\beta-1)\right)^2;\,
$$
and $P_{10,c}(L;\alpha,\beta):=\sum_{i=0}^5 p_i(\alpha,\beta)\,L^i
$, being
$$
\begin{array}{rl}
p_5(\alpha,\beta) =&\alpha \beta+1, \\
p_4(\alpha,\beta) =&5 {\alpha}^{2}{\beta}^{2}-5 {\alpha}^{2}\beta-5
\alpha {\beta}^{2}
+12 \alpha \beta-4 \alpha-4 \beta+1,\\
p_3(\alpha,\beta) =&10 {\alpha}^{3}{\beta}^{3}-20
{\alpha}^{3}{\beta}^{2}-20 {\alpha}^{ 2}{\beta}^{3}+10
{\alpha}^{3}\beta+56 {\alpha}^{2}{\beta}^{2}+10 \alpha
{\beta}^{3}-41 {\alpha}^{2}\beta-41 \alpha {\beta}^{2}\\ &+6 {
\alpha}^{2} +35 \alpha \beta+6 {\beta}^{2}-6 \alpha-6 \beta+1, \\
p_2(\alpha,\beta) =& \left( \beta-1 \right)  \left( \alpha-1 \right)
\left( 10 {\alpha} ^{3}{\beta}^{3}-20 {\alpha}^{3}{\beta}^{2}-20
{\alpha}^{2}{\beta}^{3 }+10 {\alpha}^{3}\beta+55
{\alpha}^{2}{\beta}^{2}+10 \alpha {\beta }^{3}\right.
\\
& \left.-37 {\alpha}^{2}\beta-37 \alpha {\beta}^{2}+4
{\alpha}^{2}+30
 \alpha \beta+4 {\beta}^{2}-5 \alpha-5 \beta+1 \right), \\
p_1(\alpha,\beta) =&\left( \beta \left( \alpha-1 \right) -1 \right)
^{4} \left( 5 { \alpha}^{3}{\beta}^{3}-10 {\alpha}^{3}{\beta}^{2}-10
{\alpha}^{2}{ \beta}^{3}+5 {\alpha}^{3}\beta+25
{\alpha}^{2}{\beta}^{2}+5 \alpha
 {\beta}^{3}\right.\\ &\left. -15 {\alpha}^{2}\beta-15 \alpha {\beta}^{2}+{\alpha}^{
2}+12 \alpha \beta+{\beta}^{2}-2 \alpha-2 \beta+1 \right),  \\
p_0(\alpha,\beta) =&\alpha \beta  \left( \beta-1 \right)^{5} \left(
\alpha-1 \right)^{5}.
\end{array}
$$

It is easy to check, using the expression of $5H$ given above, that
when $\alpha\beta-1\neq 0$ the value of $L$ corresponding to the
root of $P_{10,a}$ give $5H=V$, thus characterizing orbits of prime
period $5$, and so $10$ is not a prime period for the orbits in this
level set. If $\alpha \beta-1 = 0$, then $P_{10,a}$ does not depends
on $L$, and it is zero only if $\alpha=\beta=1$, which corresponds
to the well--known globally $10$ periodic case $a=b=1$.

Observe that, by using (\ref{alphatopq}), the root of $P_{10,b}$ is
$$L = \alpha(1-\beta)=\frac{p^2(1+p-q^2)}{(p+1)(q+1)},$$ which is
obviously less than $L_c$.

To  prove that the real zeros of $P_{10,c}$ are lower than $L_c$, we
proceed as in the period $6$ case, by using the change
(\ref{alphatopq}) and the translation $L=\widetilde{L}+L_c$,
obtaining that $P_{10,c}(L)=0$ if and only if the degree $5$
polynomial $ \widetilde{P}_{10,c}( \widetilde{L})$ obtained when
considering the equation $P_{10,c}(\widetilde{L}+L_c)=0$ with
parameters $p$ and $q$, vanishes. Once again it is easy to check
that this polynomial has positive coefficients when $p$ and $q$ are
positive, so its real roots are all negative, and therefore the
roots of $P_{10,c}$ are all lower than $L_c$.\qed

\subsection{Possible periods. Proof of Theorem \ref{teo2}}\label{secperiods}

In order to prove both Theorem \ref{teo2} and Corollary
\ref{reclyness} we need the following auxiliary result.

\begin{lem}\label{lemaodd}
If $a\neq b$, then the recurrence (\ref{eq})-(\ref{k=2}) do not have
periodic solutions of odd period $p\neq 1$.\end{lem}

\proof Suppose that for $u_1>0$,$u_2>0$ the sequence $\{u_n\}$ is a
periodic solution with period $p=2k+1$. Then, by using the relations
in (\ref{Fbaab}), we have that
\begin{equation}\label{primeracondicio}
\left(u_{1},u_{2}\right)=\left(u_{2k+2},u_{2k+3}\right)=F_a\left(F_{b,a}\right)^k(u_1,u_2).
\end{equation}
But on the other hand, from (\ref{Fbaab}) and
(\ref{primeracondicio}) we have:
$$
\begin{array}{rl}
F_a(u_1,u_2)&=\left(u_{2},u_{3}\right) = \left(u_{2k+3},u_{2k+4}\right)=F_b\left(F_{a,b}\right)^k(u_2,u_3)\\
   & =F_b\left(F_{a,b}\right)^k
   F_a(u_1,u_2)=\left(F_{b,a}\right)^{k+1}(u_1,u_2)\\
   &=F_b\left(F_a\left(F_{b,a}\right)^k(u_1,u_2)\right)\\
   &=F_b\left(u_1,u_2\right).
\end{array}
$$
Hence $F_a(u_1,u_2)=F_b\left(u_1,u_2\right)$, which implies that
$a=b$.\qed

\smallskip

 Also to
prove Theorem \ref{teo2}, we need the following auxiliary result
which allow us to characterize, constructively, which periods arise
when the rotation number takes values in a given interval.

\begin{lem}\label{algoperiod}(\cite[Theorem 25 and Corollary 26]{CGM07}) Consider an open interval $(c,d)$; denote by $p_1=2,p_2=3, p_3,\ldots,
p_n,\ldots$ the set of all the prime numbers, ordered following the
usual order. Also consider the following natural numbers:
\begin{itemize}
\item Let $p_{m+1}$ be the smallest prime number satisfying that $p_{m+1}>\max(3/(d-c),2),$
\item  Given any prime number $p_n,$ $1\le n\le m,$ let $s_n$ be the smallest natural number
such that $p_n^{s_n}>4/(d-c).$
\item Set  $p:=p_1^{s_1-1} p_2^{s_2-1}\cdots p_m^{s_m-1}.$
\end{itemize}
 Then, for any
$r>p$ there exists an irreducible fraction $q/r$ such that $q/r \in
(c,d).$
\end{lem}

\noindent \textsl{Proof of Theorem \ref{teo2}.} (i) If
$\sigma(a,b)\neq 2/5$  (i.e. when  $(a,b)\notin
\Gamma=\left\{\sigma(a,b)=2/5,\,a,b>0\right\}$, where $\Gamma$ is
the curve given in Corollary \ref{corolgamma}), then the value
$p_0(a,b)$ can be computed applying Lemma \ref{algoperiod} to the
interval $I(a,b)$ given in Theorem \ref{periodsinI}, as in the proof
of statement (iii) below.

If $\sigma(a,b)= 2/5$, then observe that $(a,b)=(1,1)$ is the only
point in the parameter space
$\mathcal{P}=\left\{(a,b),\,a,b>0\right\}$ such that $(a,b)\in
\left\{\Gamma\cap\{a=b\}\right\}$. Hence, as a consequence of Lemma
\ref{lemaodd}, if $(a,b)\neq (1,1)$ then  $J(a,b)=
\mathrm{Image}\left(\theta_{b,a}\left((h_c,+\infty)\right)\right)$
is a closed interval with nonempty interior (i.e. not the single
value $2/5$), so there exists some value $p_0(a,b)\in\mathbb{N}$
such that for all $p>p_0(a,b)$ there exists irreducible rational
numbers $q/p\in J(a,b)$. Otherwise $\theta_{b,a}(h)\equiv 2/5$ and
all the orbits would be $5$-periodic, a contradiction with Lemma
\ref{lemaodd}.

Observe that the computability of $p_0(a,b)$ is a generic property
in $\mathcal{P}$, since $\mathcal{P}\setminus\Gamma$ is an open and
dense subset of it.

 Statement (ii) is a direct consequence of  Proposition
\ref{unionofJ}, which implies that for each number in $(1/3,1/2)$
there exists some $a,b>0$ and some initial condition in ${\mathcal
Q}^+$  with this associated rotation number for $F_{b,a}$. In
particular, for all the irreducible rational numbers $q/p\in
(1/3,1/2)$ there are positive values of $a$ and $b$ such that
$F_{b,a}$ has continuum of periodic orbits of period $p.$

(iii) Setting $c=1/3$ and $d=1/2$,  and using the notation
introduced in Lemma \ref{algoperiod}, we have that, $m=7$; $p_1=2$
(with $s_1=5$); $p_2=3$ (with $s_2=3$); $p_3=5$; $p_4=7$; $p_5=11$;
$p_6=13$ and  $p_7=17$ (with $s_i=1$ for $i=3,\ldots,7$). From Lemma
\ref{algoperiod}, for all $p\in \N$, such that
  $p>p_0:=2^4\cdot3^3\cdot5\cdot 7\cdot 11\cdot 13\cdot
  17=12\,252\,240$ there exists an irreducible fraction
  $q/p\in(1/3,1/2)$. Hence by Proposition \ref{unionofJ} there
  exists some $a,b>0$ such that there exists a continuum of initial conditions with rotation number $\theta_{b,a}(h)=q/p$,
  thus giving rise to $p$--periodic orbits of $F_{b,a}$. Now, using a finite
algorithm we can determine which irreducible
  fractions $q/p$ with $p\le p_0$ are in $(1/2,1/3)$, resulting that there appear irreducible fractions with all  the
denominators except $2,3,4,6$ and $10$.

Finally,  the periods $2$, $3$, $4$, $6$ and $10$ do not appear as a
consequence of Propositions \ref{propoimp1} and \ref{propoimp4}.\qed

\medskip

 Notice that, as a consequence of the arguments in the
proof of Theorem \ref{teo2} (i), we have the following result:

\begin{propo}\label{nomonotonaasigma25}
Set $a,b>0$. If  $(a,b)\neq(1,1)$ and  $\sigma(a,b)=2/5$, then
$\theta_{b,a}(h)$ is not monotonic on $[h_c,+\infty)$.
\end{propo}

Now the proof of Corollary \ref{reclyness} is straightforward.

\medskip

\noindent \textsl{Proof of Corollary \ref{reclyness}.} First observe
that from Lemma \ref{lemaodd}, if $a\neq b$, then the recurrence
(\ref{eq})-(\ref{k=2}) do not have periodic solutions of odd period
$p\neq 1$.

Finally, observe that  any $p$--periodic orbit of $F_{b,a}$ with
initial conditions $(u_1,u_2)$ will give rise to a $2p$--periodic
solution of the recurrence (\ref{eq}). Now, both statements (i) and
(ii)  follow as a straightforward application of Theorem \ref{teo2}.
\qed

\end{document}